\documentclass[final]{siamltex}
\usepackage{latexsym}
\usepackage{amssymb}
\usepackage{amsfonts}
\usepackage[dvips]{graphicx}

\newtheorem{remark}[theorem]{\textit{Remark}}

\newcommand{\eps}{\varepsilon}
\newcommand{\ups}{\upsilon}

\newcommand{\R}{{\mathbb{R}}}

\begin{document}

\title{Asymptotic stability of periodic solutions for nonsmooth
differential equations with application to the nonsmooth van der Pol
oscillator}


\author{Adriana  Buic\u{a}\thanks{Department of Applied Mathematics,
        Babe\c{s}--Bolyai University, Cluj--Napoca, Romania ({\tt
        abuica@math.ubbcluj.ro}).} \and
        Jaume Llibre\thanks{Departament de Matem\`atiques,
        Universitat Aut\`onoma de Barcelona, 08193
        Bellaterra, Barcelona, Spain ({\tt jllibre@mat.uab.cat}).} \and
        Oleg Makarenkov\thanks{Research Institute of Mathematics, Voronezh State
        University, Voronezh, Russia ({\tt omakarenkov@math.vsu.ru}).}}

\renewcommand{\labelenumi}{(\theenumi)}
\renewcommand{\theenumi}{\roman{enumi}}

\maketitle

\begin{abstract}
In this paper we study the existence, uniqueness  and asymptotic
stability of the periodic solutions for the  Lipschitz system $
\dot x =\eps g(t,x,\eps).$ Classical hypotheses in the periodic
case of second Bogolyubov's theorem imply our ones. By means of
the results established we construct, for small $\eps,$ the curves
of dependence of the amplitude of asymptotically stable
$2\pi$--periodic solutions of the nonsmooth van der Pol oscillator
$ \ddot u+\eps \left(|u|-1\right)\dot u+(1+a\eps)u=\eps\lambda\sin
t, $ on the detuning parameter $a$ and the amplitude of the
perturbation $\lambda.$ After, we compare the resonance curves
obtained, with the resonance curves of the classical van der Pol
oscillator $ \ddot u+\eps \left(u^2-1\right)\dot
u+(1+a\eps)u=\eps\lambda\sin t,$ which were first constructed by
Andronov and Witt.
\end{abstract}



\pagestyle{myheadings}

\thispagestyle{plain}

\markboth{A. BUIC\u{A}, J. LLIBRE AND O. MAKARENKOV}{PERIODIC
SOLUTIONS FOR NONSMOOTH DIFFERENTIAL SYSTEMS}

\section{Introduction}

\newcounter{items}
\setcounter{items}{0}

In the present paper we study the existence, uniqueness and asymptotic stability of the $T$--periodic solutions for the system
\begin{equation}\label{ps}
\dot x =\eps g(t,x,\eps),
\end{equation}
where $\eps>0$ is a small parameter and the function $g\in C^0(\mathbb{R} \times\mathbb{R}^k\times[0,1],\mathbb{R}^k)$ is $T$--periodic in the first variable and locally  Lipschitz with respect to the second one. As usual a key role will be played by the averaging
function
\begin{equation} \label{averfunction}
{g}_0(v)=\int\limits_0^T g(\tau,v,0)d\tau,
\end{equation}
and we shall look for those periodic solutions  that starts near some $v_0\in {g}_0^{-1}(0)$.

In the case that $g$ is of class $C^1$, we remind the periodic case of the second Bogolyubov's theorem (\cite{bog}, Ch.~1, \S~5, Theorem~II) which represents a part of the averaging principle: {\it $\det \, (g_0)'(v_0)\neq 0$
assures the existence and uniqueness, for $\eps>0$ small, of a $T$--periodic
solution of system (\ref{ps}) in a neighborhood of $v_0,$ while the fact that all the eigenvalues of the Jacobian matrix $(g_0)'(v_0)$ have negative real part, provides also its asymptotic stability.}
 This theorem has a long history and it includes results by Fatou \cite{fat}, Mandelstam--Papaleksi  \cite{mand} and Krylov--Bogolyubov \cite[\S~2]{kry}.

\begin{figure}[h]
\begin{center}
\includegraphics[scale=0.2]{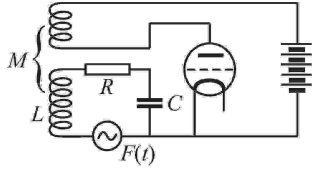}\end{center}
\caption{\footnotesize Circuit scheme for the classical lamp
oscillator (see Andronov-Vitt-Khaikin \cite{andb}, Ch.VIII, \S2,
Fig.~348, Malkin \cite{mal}, Ch.I, \S5, Fig.~1, Nayfeh-Mook
\cite{nay}, \S3.1.7, Fig.~3-5).} \label{thom}
\end{figure}

Second Bogolubov's theorem gave a theoretical justification of
resonance phenomenons in many real physical systems. The most
significant example is the classical lamp oscillator  whose scheme
is drawn at Fig.~\ref{thom} and whose current $u$ is described by
the following second order differential equation
\begin{equation}\label{two}
  \ddot u+\frac{1}{LC}\left(RC-Mi'(u)\right)\dot u+\omega^2 u=\frac{1}{LC}F(t),
\end{equation}
where $R=\eps R_0,$ $M=\eps M_0,$ $\omega^2=1+\eps b,$
$F(t)=\eps\lambda\sin t,$ $\eps>0$ is assumed to be small and the
lamp characteristic is drawn
\begin{figure}[h]
\begin{center}
\includegraphics[scale=0.6]{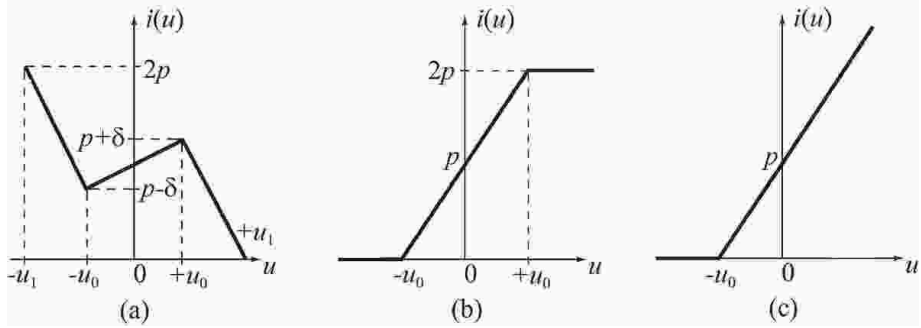}\end{center}
\caption{\footnotesize Characteristics of the lamp of the circuit
of Fig.~\ref{thom}. (a) - lamp in a harsh regime (see
Andronov-Vitt-Khaikin \cite{andb}, Ch.IV, \S7, Fig.212b, Malkin
\cite{mal}, Ch.I, \S5, comments for Eq.~5.3-5.4); (b) - lamp with
saturation (see Andronov-Vitt-Khaikin \cite{andb}, Ch.VIII, \S3,
Fig.~364); (c) - lamp without saturation (see
Andronov-Vitt-Khaikin \cite{andb}, Ch.IX, \S7, Fig.~482)}
\label{char}
\end{figure}
at Fig~\ref{char}a. The analysis of bifurcation of periodic
solutions in this system is performed in almost every book on
nonlinear oscillations (see Andronov-Vitt-Khaikin \cite{andb},
Ch.VIII, \S2, Malkin \cite{mal}, Ch.I, \S5, Nayfeh-Mook
\cite{nay}, \S3.1.7) but with approximation
$i(u)=i_{(a)}(u)=S_0+S_1 u-\frac{1}{3}S_3 u^3$ (leading to the
classical van der Pol equation). Therefore, it is natural to look
for a technique that permits avoiding this approximation and
dealing with the original shape of the lamp characteristic drawn
at Fig~\ref{char}a that expects to give more accurate
correspondence between theoretical and experimental results.
Moreover, a wide class of physical systems is modelled by circuit
Fig~\ref{thom} whose lamp either has or has not a saturation that
leads to the characteristic drawn at Figures \ref{char}b and
\ref{char}c respectively. Though the unforced  equation
(\ref{two}) (i.e. for $F=0$) with $i$ described by
Fig.~\ref{char}b and Fig.~\ref{char}c is well studied (see
Andronov-Vitt-Khaikin \cite{andb}, Ch.VIII, \S3 and Ch.IX, \S7
respectively), the question about resonances in these equations
when $F\not=0$ (e.g. $F(t)=\eps\lambda\sin t$) is  still open.
With regard to equation (\ref{two})
 with lamp
characteristic given by Fig.~\ref{char}a, \ref{char}b or
\ref{char}c we finally note that Levinson's change of variables
(see \cite{lev}, pass from Eq.~2.0 to Eq.~2.1) allows to rewrite
equation (\ref{two}) as the following system
\begin{eqnarray*}
& & \dot z_1=z_2-\frac{z_1 RC +M(i(z_1)-i(0))}{LC}=z_2-\eps\frac{z_1  R_0 C + M_0(i(z_1)-i(0))}{LC}\\
 & &  \dot z_2=-\omega z_1+F(t)=-z_1+\eps\left(b z_1+\frac{1}{LC}\lambda\sin t\right),
\end{eqnarray*} whose solution $(z_1,z_2)$ gives a solution $u=z_1$ to
(\ref{two}). Then the change of variables $
  \left(\begin{array}{cc}
    z_1(t) \\
    z_2(t)
    \end{array}\right)=\left(\begin{array}{cc}
    \cos t & \sin t \\
    -\sin t & \cos t
    \end{array}\right)\left(\begin{array}{cc}
    x_1(t) \\
    x_2(t)
    \end{array}\right)
$ brings this system to the form (\ref{ps}) with Lipschitzian (in
the second variable) $g.$ Therefore the goal of the paper is to
generalize second Bogolubov's theorem for the case when $g$ in
(\ref{ps}) is Lipschitzian.

\begin{figure}[h]
\begin{center}
\includegraphics[scale=0.4]{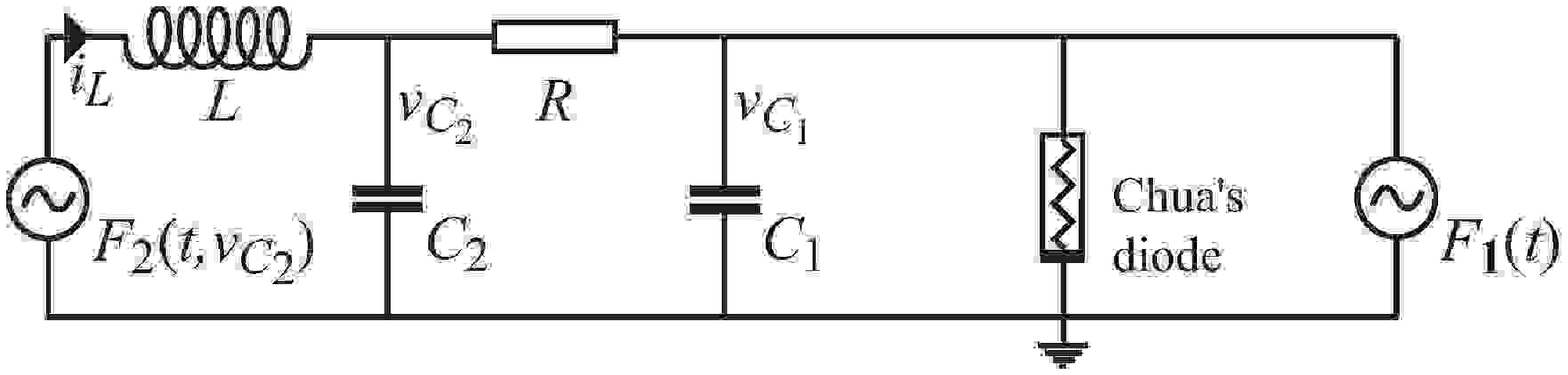}\end{center}
\caption{\footnotesize Driven Chua's circuit (see \cite{good2},
\cite{chua1}, \cite{good3}, \cite{good4}, \cite{good1}).}
\label{chu}
\end{figure}

\begin{figure}[h]
\begin{center}
\includegraphics[scale=0.23]{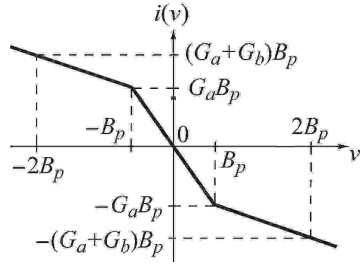}\end{center}
\caption{\footnotesize Nonlinear characteristic of the Chua's
diode of the circuit drawn at Fig.~\ref{chu} given by $i(v)=G_b
v+(1/2)(G_a-G_b)\left(\left|v+B_p\right|-\left|v-B_p\right|\right),$
 where $G_a,G_b,B_p\in\mathbb{R}$ are some constants
 depending on the properties of the Chua's diode (see Chua \cite{ch})} \label{charchu}
\end{figure}

Another motivation of this paper comes from modern electrical
engineering where driven (or forced) Chua's circuit drawn at
Fig.~\ref{chu} is a subject of an enormous number of papers.
Circuit at Fig.~\ref{chu} is described by the following
three-dimensional system
\begin{eqnarray}
 C_1
 \frac{dv_{C_1}}{dt}&=&\frac{v_{C_2}-v_{C_1}}{R}-i(v_{C_1})+F_1(t),\nonumber\\
 C_2\frac{dv_{C_2}}{dt}&=&\frac{v_{C_1}-v_{C_2}}{R}+i_L,\label{chua}\\
 L \frac{di_L}{dt}&=&-v_{C_2}+F_2(t,v_{C_2})\nonumber
\end{eqnarray}
where $i(v)$ is the characteristic of the Chua's diode whose shape
drawn at Fig.~\ref{charchu} is piecewise linear.
Many numerical simulations have been suggested around dynamics of
(\ref{chua}) in the recent literature, see \cite{good1},
\cite{good3} for
 $F_1\not=0$ and $F_2\not=0,$  \cite{good2}, \cite{good4} for $F_1=0$ and
periodic $F_2$, \cite{chua1} where both $F_1$ and $F_2$ are
periodic. Generalization of the second Bogolubov's theorem for
equations (\ref{ps}) with Lipschitzian right hand part will allow
for the first time theoretical detection of asymptotically stable
periodic solutions in certain configurations of the driven Chua's
circuit (\ref{chua}) provided that $C_1$ is large enough. This
promises to forestall some numerical simulations (e.g. to work out
interesting parameters of the driven Chua's circuit) giving a
significant impact for further experiments.

\begin{figure}[h]
\begin{center}
\includegraphics[scale=0.4]{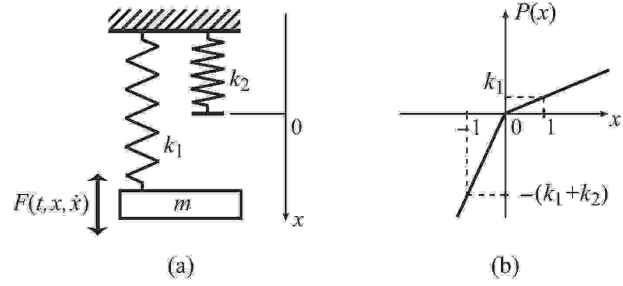}\end{center}
\caption{\footnotesize  A prototypic device (a) where a driven
mass is attached to a immovable beam via a spring with piecewise
linear stiffness (b), see e.g. \cite{bov},   \cite{kru} (Ch.I,
p.16 and Ch.IV, p.100) and \cite{ieee}.} \label{spr}
\end{figure}

For a large extent the phenomenal interest in  generalizing of the
second Bogolubov's theorem for the Lipschitzian case comes from
mechanics, where systems with piecewise linear stiffness describe
various oscillating processes. A prototypic process of this type
is exhibited by the device drawn at Fig.~\ref{spr}a where a forced
mass is attached to a spring whose stiffness changes from $k_1$ to
$k_1+k_2$ when the mass coordinate crosses $0$ in the negative
direction. This device is governed by the following second order
differential equation
\begin{equation}\label{spring}
  m\ddot x+P(x)=F(t,x,\dot x),
\end{equation}
where piecewise linear stiffness $P$ is drawn at Fig.~\ref{spr}b.
Depending on a particular configuration of the device of
Fig.~\ref{spr}a various terms can stay for $F$ in (\ref{spring}).
It is $F(t,x,\dot x)=-f(x)\dot x+M\cos\omega t$ with piecewise
constant $f$ for shock-absorber and jigging conveyor (see
\cite{kru}, Ch.I, p.16 and Ch.IV, p.100 where original second
Bogolubov's theorem is employed without justification). Levinson's
change of variables (see transformation of Eq.~(\ref{two}) above)
allows to rewrite (\ref{spring}) as a Lipschitz system. It takes
simpler form $F(t,x,\dot x)=-c\dot x+M Q(t)$ for an impact
resonator  and $F(t,x,\dot x)=-c\dot x+M\sin\omega t$ for a
cracked-body model (see \cite{ieee} and \cite{bov}, where
numerical experiments are performed solely). In each of these
situations equation (\ref{spring}) can be rewritten as system
(\ref{ps}) with Lipschitzian $g$ provided that the constant $k_2$
and the amplitude of the force $F$ are sufficiently small.
Therefore the related generalization of the second Bogolubov's
theorem promises to justify or explain the resonances appeared in
the cited results. We note that the recent report by Los Alamos
National Laboratory \cite{alamos} describes increasing interest in
a specific form of the model of Fig.~\ref{spr}a called
cracked-body model and, particularly, in suspension bridges models
that is why we discuss the contribution of the present paper to
the latter model in details.

\begin{figure}[h]
\begin{center}
\includegraphics[scale=0.5]{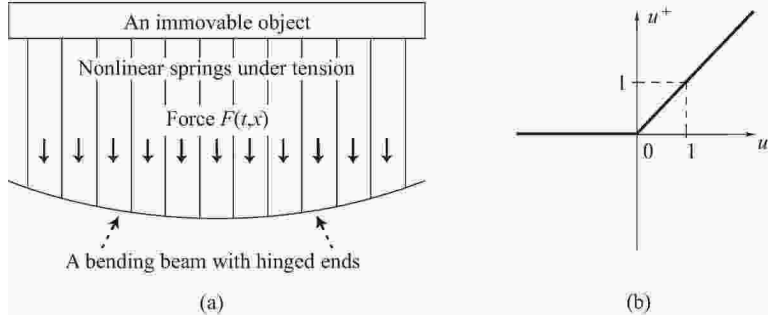}\end{center}
\caption{\footnotesize (a) -- The first idealization of the
suspension bridge: the beam bending under its own weight is
supported by the nonlinear cables (see Lazer-McKenna \cite{laz},
Fig.~2); (b) -- characteristic of stiffness of nonlinear springs.}
\label{bri}
\end{figure}
The first idealization of a one-dimensional suspended bridge is
drawn at Fig.~\ref{bri}a. It is represented (see \cite{zamp},
\cite{laz}) by the beam bending under its own weight and supported
by cables whose restoring force due to elasticity is proportional
to $u^+$ (see Fig.~\ref{bri}b), where $u=u(t,x)$ is the
displacement at a point at distance $x$ from one end of the bridge
at time $t$ and $u$ measured in the downward direction. Looking
for $u$ of the form $u(t,x)=z(t)\sin(\pi x/L)$ and considering
$F(x,t)=h(t)\sin(\pi x/L)$ we arrive (see \cite{zamp}) to the
following form of equation (\ref{spring})
\begin{equation}\label{bridge}
  m\ddot z+\delta\dot z+c(\pi/L)^4z+dz^+=mg+h(t),
\end{equation}
where the constant $m>0$ is the mass per unit length, $\delta>0$ is
a small viscous damping coefficient, $c>0$ measured the flexibility
or stiffness of the bridge, $L>0$ is the length of the bridge, $d>0$
represents the stiffness of nonlinear springs and $h$ is a
continuous $T$-periodic force modelling wind, marching troops or
cattle (see \cite{report} for details). Considering $c>0$ and $d>0$
fixed and assuming that either $c>0$ and $h(t)$ are sufficiently
small, or $c>0$ is fixed and $h(t)$ is sufficiently large, or $c>0$
is sufficiently small and $h(t)$ fixed, Giover, Lazer, McKenna,
Fabry (see \cite{zamp}, \cite{tams}, \cite{laz}, \cite{fabry})
proved various theorems on location of asymptotically stable
$T$-periodic solutions in (\ref{bridge}). The question what happens
with these solutions when $d>0$ appears to be the same magnitude
small as $\delta>0$ and $h(t)$, is open for a while and can be
resolved by means of the generalization of the second Bogolubov's
theorem we propose. Lazer and McKenna proved in \cite{tams} that the
Poincar\'e map for (\ref{bridge}) is differentiable, but we note
that it is not an enough argument to apply the original second
Bogolubov's theorem since it requires the differentiability of
functions participating in (\ref{bridge}) as well.

At the end of the applications review we note that system
(\ref{chua}) describing the Chua's circuit (Fig.~\ref{chu})
appeared in the recent time to govern mechanical systems with
so-called "negative slope" (see Awrejcewicz \cite{awr}, \S8.2.2).
Thus, analogous to Chua's circuit applications of the result of
this paper are also possible to these mechanical systems.

It was Mitropol'skii  who first noticed that various applications
require generalization of the second Bogolubov's theorem for
Lipschitz right hand parts. Assuming that $g$ is Lipschitz,
$g_0\in C^3(\R^k,\R^k)$ and that all the eigenvalues of the matrix
$(g_0)'(v_0)$ have negative real part Mitropol'skii developed the
second Bogolyubov's Theorem proving the existence and uniqueness
of a $T$--periodic solution of system (\ref{ps}) in a neighborhood
of $v_0.$ There was a great progress weakening the assumptions of
Mitropol'skii in his existence result (see Samoylenko \cite{sam}
and Mawhin \cite{maw1}), but not of his uniqueness result.
Moreover, the asymptotic stability  conclusion of the second
Bogolyubov's Theorem remained to be not generalized for
Mitropol'skii's settings (namely, when $g$ is Lipschitz) for a
long time. It has been done recently by Buic\u{a}--Daniilidis in
\cite{adr} for a class of functions $v\mapsto g(t,v,0)$
differentiable at $v_0$ for almost any $t\in[0,T],$ but it is
assumed in \cite{adr} that the eigenvectors of the matrix
$({g}_0)'(v_0)$ are orthogonal.

In the next section of the paper assuming that $g$ is piecewise
differentiable in the second variable we show in Theorem~\ref{th1}
that Mitropol'skii's conditions imply not only uniqueness, but
also asymptotic stability of a $T$--periodic solution of system
(\ref{ps}) in a neighborhood of $v_0.$ In other words we show that
Bogolyubov's theorem formulated above is valid when $g$ is  not
necessary $C^1$. Theorem~\ref{th1} follows from our even more
general Theorem~\ref{th0} whose hypotheses do not use any
differentiability neither of $g$ nor of $g_0$.
 In Section~3 we illustrate our result constructing resonance curves of
 nonsmooth van der Pol oscillator \cite{hog}. This application has
 been chosen since it allows to compare the issues of our Theorem~\ref{th1} with the
 classical results \cite{andr} and \cite{andr1}  by
Andronov and Witt obtained for original van der Pol oscillator.

\section{Main results} Throughout the paper $\Omega\subset\R^k$ is some
open set. For any $\delta>0$ we denote $B_\delta(v_0)=\left\{ v\in
\R^k~:~\|v-v_0\|\leq \delta \right\}$.
We have the following main result.

\begin{theorem}\label{th0}
Let $g\in C^0(\mathbb{R}\times \Omega \times [0,1],\mathbb{R}^k)$ and $v_0\in \Omega$. Assume the following four conditions.
\begin{enumerate}
\setcounter{enumi}{\value{items}}
\item\label{lip} For some $L>0$ we have that
$\left\|g(t,v_1,\eps)-g(t,v_2,\eps)\right\|\le
L\left\|v_1-v_2\right\|$
for any $t\in[0,T],\ v_1,v_2\in \Omega,$
$\eps\in[0,1].$
\setcounter{items}{\value{enumi}}
\item\label{newunif} For any $\gamma>0$ there exists $\delta>0$ such that
\[
\begin{array}{l}
\hskip-1cm\left\|\int_0^T g(\tau,v_1+u(\tau),\eps)d\tau-\int_0^T
g(\tau,v_2+u(\tau),\eps)d\tau\right.\\
\left.-\int_0^T g(\tau,v_1,0)d\tau+\int_0^T
g(\tau,v_2,0)d\tau\right\|\le\gamma\|v_1-v_2\|
\end{array}
\]
for any $u\in C^0([0,T],\mathbb{R}^k),$ $\|u\|\le\delta,$
$v_1,v_2\in B_\delta(v_0)$ and $\eps\in[0,\delta].$
\setcounter{items}{\value{enumi}}
\item\label{def}Let $g_0$ be the averaging function given by
(\ref{averfunction}) and consider that ${g}_0(v_0)=0.$
\setcounter{items}{\value{enumi}}
\item\label{conde}
There exist $q\in[0,1),$ $\alpha,\delta_0>0$ and a norm
$\|\cdot\|_0$ on $\R^k$  such that $\left\|v_1+ \alpha
{g}_0(v_1)\right.$ $\left. -v_2- \alpha {g}_0(v_2)\right\|_0 \leq
q \|v_1-v_2\|_0$ for any $v_1,v_2\in  B_{\delta_0}(v_0).$
\setcounter{items}{\value{enumi}}
\end{enumerate}
Then there exists $\delta_1>0$ such that for every
$\eps\in(0,\delta_1]$ system (\ref{ps}) has exactly one
$T$--periodic solution $x_\eps$ with $x_\eps(0)\in
B_{\delta_1}(v_0).$  Moreover the solution $x_\eps$ is
asymptotically stable and $x_\eps(0)\to v_0$ as $\eps\to 0.$
\end{theorem}

When solution $x(\cdot,v,\eps)$ of system (\ref{ps}) with initial condition
$x(0,v,\eps)=v$ is well defined on $[0,T]$ for any $v\in
B_{\delta_0}(v_0)$, the map $v\mapsto x(T,v,\eps)$ is well defined
and it is said to be the {\it Poincar\'e map} of system (\ref{ps}).
The proof of existence, uniqueness and stability of the
$T$--periodic solutions of system (\ref{ps}) in Theorem~\ref{th0} reduces
to the study of corresponding properties of the fixed points of this
map.

In order to prove Theorem \ref{th0} we observe from (\ref{ps}) that
$x(T,v,\eps)$ can be represented as
\[
x(T,v,\eps)=v+\eps g_\eps(v),\ \ {\rm where }\ \
g_\eps(v)=\int\limits_0^T g(\tau,x(\tau,v,\eps),\eps)d\tau,
\]
and we use the following result which claims that properties (\ref{lip})
and (\ref{newunif}) are also applied to $g_\eps$ in a suitable sense.

\begin{lemma}\label{lem3}
Let $g\in C^0(\mathbb{R}\times\Omega\times[0,1],\mathbb{R}^k)$ and
$\delta_0>0$ be such that $B_{\delta_0}(v_0)\subset\Omega.$ If
(\ref{lip}) is satisfied then there exist $\delta\in[0,\delta_0]$
and $L_1>0$ such that the map $(v,\eps)\mapsto g_\eps(v)$ is well
defined and continuous on $B_{\delta_0}(v_0)\times [0,\delta]$ and
\[
\|g_\eps(v_1)-g_\eps(v_2)\|\le L_1\|v_1-v_2\|\quad{ for\ any\
}\eps\in[0,\delta],\ v_1,v_2\in B_{\delta_0}(v_0).
\]
If both (\ref{lip}) and (\ref{newunif}) are satisfied then for any
$\gamma>0$ there exists $\delta\in[0,\delta_0]$ such that
\[
\| g_\eps(v_1)-g_0(v_1)-g_\eps(v_2)+g_0(v_2)\| \le \gamma \| v_1-v_2
\| \] for any $v_1,v_2\in B_\delta(v_0)$ and $\eps\in[0,\delta].$
\end{lemma}

\begin{proof}
Using the continuity of the solution of a differential system with
respect to the initial data and the parameter (see \cite{pon},
Ch.~4, \S~23, statements G and D), we obtain the existence of
$\eps_0>0$ such that $x(t,v,\eps)\in\Omega$ for any $t\in[0,T],$
$v\in B_{\delta_0}(v_0)$ and $\eps\in[0,\eps_0].$ Using the
Gronwall--Bellman Lemma \cite[Ch.~II, \S~11]{dem} from the
representation $x(t,v,\eps)=v+\eps\int\limits_0^t
g(\tau,x(\tau,v,\eps),\eps)d\tau$ and the property (\ref{lip}) we obtain
$ \|x(t,v_1,\eps)-x(t,v_2,\eps)\|\leq e^{\eps LT} \|v_1-v_2\|$ for
all $t\in[0,T],$ $v_1, \, v_2 \in B_{\delta_0}(v_0)$ and
$\eps\in[0,\eps_0].$ Therefore $y(t,v,\eps)=\int_0^t
g(\tau,x(\tau,v,\eps),\eps)d\tau$ satisfies the following property
\begin{equation}
\label{lip-y} \|y(t,v_1,\eps)-y(t,v_2,\eps)\|\leq L_1 \|v_1-v_2\|
\end{equation}
$\mbox{for all } t\in[0,T],\ v_1,\ v_2 \in B_{\delta_0}(v_0),\
\eps\in[0,\eps_0]$ and $L_1= LTe^{\eps_0 LT}.$
Since $g_\eps(v)=y(T,v,\eps)$ the first part of the
lemma has been proven.

Taking into account that
$x(t,v,\eps)=v+\eps y(t,v,\eps)$ we have
\begin{equation}\label{SU}
\ \ \ \ \ y(T,v_1,\eps)- y(T,v_1,0)-y(T,v_2,
\eps)+y(T,v_2,0)=I_1(v_1,v_2,\eps)+I_2(v_1,v_2,\eps)
\end{equation}
 where
\begin{eqnarray*}
I_1(v_1,v_2,\eps)&=& \int_0^T [g(\tau,v_2+\eps
y(\tau,v_1,\eps),\eps)- g(\tau,v_2+\eps y(\tau,v_2,\eps),
\eps)]d\tau
\\
I_2(v_1,v_2,\eps)&=& \int_0^T [(g(\tau,v_1+\eps
y(\tau,v_1,\eps),\eps)- g(\tau,v_2+\eps
y(\tau,v_1,\eps),\eps))]d\tau\\
& &\qquad\qquad\qquad-\int_0^T{(g(\tau,v_1,0)-
g(\tau,v_2,0))}d\tau.
\end{eqnarray*}
Since $(t,\ups,\eps)\mapsto y(t,\ups,\eps)$  is bounded on
$[0,T]\times B_{\delta_0}(v_0) \times [0,\eps_0]$, we have that
$\eps y(t,\ups,\eps)\to 0$ as $\eps\to 0$ uniformly with respect to
$t\in[0,T]$ and $v\in B_{\delta_0}(v_0).$ Decreasing $\eps_0>0,$ if
necessary, we get that $v_2+\eps y(t,v_1,\eps)\in\Omega$ for
any $t\in[0,T],$ $v_1,v_2\in B_{\delta_0}(v_0),$
$\eps\in[0,\eps_0].$ By assumption (\ref{lip}) and relation
(\ref{lip-y}) we obtain that $\|I_1(v_1,v_2,\eps)\|\le T\cdot\eps
LL_1\|v_1-v_2\|$ for all $\eps\in [0,\eps_0]$, $v_1,v_2\in
B_{\delta_0}(v_0).$

We fix $\gamma>0$ and take $\delta>0$ given by (\ref{newunif}).
Without loss of generality we can consider that $\delta\le
\min\{\delta_0,\, \eps_0,\, \gamma/(2TLL_1)\}.$ Therefore assumption
(\ref{newunif}) implies
that $\|I_2(v_1,v_2,\eps)\|\le (\gamma/2)\|v_1-v_2\|$ for any
$\eps\in[0,\delta],$ $v_1,v_2\in B_\delta(v_0)$. Substituting the
obtained estimations for $I_1$ and $I_2$ into (\ref{SU}) we have
$\|y(T,v_1,\eps)- y(T,v_1,0)-y(T,v_2, \eps)+y(T,v_2,0)\|\le (\eps
TLL_1+\gamma/2)\|v_1-v_2\|\le\gamma\|v_1-v_2\|$ for any
$\eps\in[0,\delta],$ $v_1,v_2\in B_\delta(v_0)$. Hence the proof is
complete.
\end{proof}

\begin{lemma} \label{bOleg}
Let $g_0:\Omega\to\R^k$  satisfy assumption (\ref{conde}) with some $q\in(0,1),$ $\alpha,\delta_0>0$ and a norm $\|\cdot\|_0$ on $\R^k.$ Then
$\left\|v_1+\eps {g}_0(v_1)-v_2-\eps {g}_0(v_2)\right\|_0\le
\left(1-\eps (1-q)/ \alpha \right)\|v_1-v_2\|_0$ \ for any $v_1,v_2\in
B_{\delta_0}(v_0)$  and any $\eps \in [0,\alpha].$
\end{lemma}

\begin{proof}
Indeed, the representation $v+\eps g_0(v)=(1-\eps/ \alpha)v + \eps /\alpha \,(v + \alpha g_0(v))$ implies that the Lipschitz constant of the function
$I+\eps g_0$  with respect to the norm $\| \cdot \| _0$ is
$(1-\eps/ \alpha) + \eps /\alpha \, q= 1-\eps (1-q)/ \alpha$.
\end{proof}

{\it Proof of Theorem~\ref{th0}.} By Lemma~\ref{lem3} we have that
there exists $\delta_1\in[0,\delta_0]$  such that
\begin{equation}\label{WU}
  \|g_\eps(v_1)-g_0(v_1)-g_\eps(v_2)+g_0(v_2)\|_0\le((1-q)/(2\alpha))\|v_1-v_2\|_0
\end{equation}
for any $\eps\in[0,\delta_1],$ $v_1,v_2\in B_{\delta_1}(v_0).$ First
we prove that there exists $\eps_1\in[0,\delta_1]$ such that for
every $\eps \in [0,\eps_1]$ there exists $v_\eps\in
B_{\delta_1}(v_0)$ such that $x(\cdot,v_\eps,\eps)$ is a
$T$--periodic solution of (\ref{ps}) by showing that there exists
$v_\eps$  such that $x(T,v_\eps,\eps)=v_\eps$. Using (\ref{def}) and
(\ref{conde}) we have
$$
  \|v+\alpha g_0(v)-v_0\|_0\le q\|v-v_0\|_0\quad{ \rm for\ any\ }v\in
  B_{\delta_1}(v_0).
$$
Therefore we have that the map $I+\alpha g_0$
 maps $B_{\delta_1}(v_0)$
into itself. From Lemma~\ref{lem3} we have that there exists
$\eps_0>0$ such that the map $(v,\eps)\mapsto g_\eps(v)$ is well
defined and continuous on $B_{\delta_1}(v_0)\times [0,\eps_0]$. We
deduce that there exists $\eps_1>0$ sufficiently small such that,
for every $\eps\in[0,\eps_1]$, the map $I+\alpha g_\eps$ maps
$B_{\delta_1}(v_0)$ into itself as well. Therefore, by the Brouwer
Theorem (see, for example, \cite[Theorem~3.1]{kraop}) we have that
$B_{\delta_1}(v_0)$ contains at least one fixed point of the map
$I+\alpha g_\eps$ for any $\eps\in[0,\eps_1].$ Denote this fixed
point by $v_\eps$. Then we have $g_\eps(v_\eps)=0$ and
$x(T,v_\eps,\eps)=v_\eps$ for any $\eps\in[0,\eps_1]$.

Now we prove that  $x(\cdot,v_\eps,\eps)$ is the only $T$--periodic
solution of (\ref{ps}) originating near $v_0$ and that, moreover, it
is asymptotically stable. Knowing that $x(T,v,\eps)=v+\eps
g_\eps(v)$ we write the following identity
\begin{equation} \label{Oleg} x(T,v,\eps)=v+\eps g_0(v)+\eps\left(
g_\eps(v)-g_0(v)\right). \end{equation}

Using Lemma~\ref{bOleg} we have from (\ref{WU}) and (\ref{Oleg})
that
\[
\begin{array}{rl}
  \left\|x(T,v_1,\eps)-x(T,v_2,\eps)\right\|_0\le &
  (1-\eps(1-q)/\alpha+\eps (1-q)/(2\alpha))\|v_1-v_2\|_0 \\
  =& (1-\eps(1-q)/(2\alpha))\|v_1-v_2\|_0,
\end{array}
\]
for all $v_1,v_2\in B_{\delta_1}(v_0)$ and $\eps\in[0,\delta_1]$. We
proved before that there exists $\eps_1>0$ that, for every
$\eps\in[0,\eps_1]$ there exists $v_\eps\in B_{\delta_1}(v_0)$ such
that $x(\cdot,v_\eps,\eps)$ is a $T$--periodic solution of
(\ref{ps}).
  Since $\eps(1-q)/(2\alpha)>0$ and  $\eps_1\leq \delta_1$ the last inequality implies
   that for each $\eps\in[0,\delta_1]$,  the $T$--periodic solution
$x(\cdot,v_\eps,\eps)$ is the only $T$--periodic solution of
(\ref{ps}) in $B_{\delta_1}(v_0)$ and, moreover (see
\cite[Lemma~9.2]{kraop}) it is asymptotically stable. $\Box$

\begin{remark}\label{ad1}
We note that a similar result close to Theorem~\ref{th0} is obtained by
Buic\u{a} and Daniilidis (see \cite{adr}, Theorem~3.5). But instead
of the assumption (\ref{conde}) with fixed $\alpha>0$ it is assumed
to be satisfied for any $\alpha>0$ sufficiently small. Although,
Lemma~\ref{bOleg} now implies that it is the same to assume
(\ref{conde}) for only one $\alpha>0$ and, respectively, for all
$\alpha>0$ sufficiently small. The advantage of our
Theorem~\ref{th0} is that it does not require differentiability of
$g(t,\cdot,\eps)$ at any point, while \cite{adr} needs it at $v_0.$
See also Remark~\ref{ad2}.
\end{remark}

In general it is not easy to check assumptions (\ref{newunif}) and
(\ref{conde}) in the applications of Theorem \ref{th0}. Thus we
give also the following theorem based on Theorem~\ref{th0} which
assumes certain type of piecewise differentiability instead of
(\ref{newunif}) and deals with properties of the matrix
$(g_0)'(v_0)$ instead of the Lipschitz constant of $g_0.$

For any set $M\subset [0,T]$ measurable
in the sense of Lebesgue we denote by ${\rm mes}(M)$ the Lebesgue
measure of $M$ (see \cite{kol}, Ch.~V, \S~3).

\begin{theorem}\label{th1}
Let $g\in C^0(\mathbb{R}\times\Omega\times[0,1],\mathbb{R}^k)$ satisfy
(\ref{lip}). Let $g_0$ be the averaging function given by
(\ref{averfunction}) and consider $v_0\in\Omega$ such that
${g}_0(v_0)=0.$ 
Assume that
\begin{enumerate}
\setcounter{enumi}{\value{items}}
\item\label{newdef} given
 any $\widetilde{\gamma}>0$  there exist $\widetilde{\delta}>0$ and $M\subset [0,T]$
measurable in the sense of Lebesgue with ${\rm
mes}(M)<\widetilde{\gamma}$ such that for every $v\in
B_{\widetilde{\delta}}(v_0),$  $t\in[0,T]\setminus M$ and
$\eps\in[0,\widetilde{\delta}]$ we have that $g(t,\cdot,\eps)$ is
differentiable at $v$ and
$\|g'_v(t,v,\eps)-g'_v(t,v_0,0)\|\le\widetilde{\gamma}$.

 \setcounter{items}{\value{enumi}}
\end{enumerate}
Finally assume that
\begin{enumerate}
\setcounter{enumi}{\value{items}}
\item\label{NE0o}
$g_0$ is continuously differentiable in a neighborhood of $v_0$ and
the real parts of all the eigenvalues of $({g}_0)'(v_0)$ are
negative. \setcounter{items}{\value{enumi}}
\end{enumerate}
Then there exists $\delta_1>0$ such that for every
$\eps\in(0,\delta_1]$,  system (\ref{ps}) has exactly one
$T$--periodic solution $x_\eps$ with $x_\eps(0)\in
B_{\delta_1}(v_0).$ Moreover the solution $x_\eps$ is asymptotically stable and $x_\eps(0)\to v_0$ as $\eps\to 0.$
\end{theorem}

For proving Theorem \ref{th1} we need two preliminary lemmas.

\begin{lemma}\label{lemma1}
Let $g\in C^0(\mathbb{R}\times\Omega\times[0,1],\mathbb{R}^k)$ satisfying
(\ref{lip}). If  (\ref{newdef}) holds then (\ref{newunif}) is satisfied.
\end{lemma}

\begin{proof}
Let $\gamma>0$ be an arbitrary number. We show that
(\ref{newunif}) holds with $\delta=\widetilde{\delta}/2,$ where
$\widetilde{\delta}$ is given by (\ref{newdef}) applied with
$\widetilde{\gamma}=\min\{\gamma/(4L),\gamma/(4T)\}.$ We consider
also $M\subset [0,T]$ given by (\ref{newdef}) applied with the
same value of  $\widetilde{\gamma}$.

Let $u\in C^0([0,T],\mathbb{R}^k),$ $\|u\|\le\delta$ and
$F(v)=\int_0^T g(\tau,v+u(\tau),\eps)d\tau-\int_0^T
g(\tau,v,0)d\tau.$ 
Let $v_1,v_2\in B_\delta(v_0)$  and $\varepsilon\in[0,\delta]$. 
We have ${F}(v)={F}_1(v)+{F}_2(v),$ where ${F}_1(v)=\int_M
(g(\tau,v+u(\tau),\eps)-g(\tau,v,0))d\tau$ and
${F}_2(v)=\int_{[0,T]\setminus M }
(g(\tau,v+u(\tau),\eps)-g(\tau,v,0))d\tau$.
 By (\ref{lip}) we have that
$\|{F}_1(v_1)-{F}_1(v_2)\|\le 2L\cdot {\rm mes}(M)\|v_1-v_2\|<
2L\widetilde{\gamma}\|v_1-v_2\|\le(\gamma/2)\|v_1-v_2\|.$ On the
other hand, using (\ref{newdef}), we will prove that a similar
relation holds for $F_2$. In order to do this, we denote
$h(\tau,v)=g(\tau,v+u(\tau),\eps)-g(\tau,v,0)$. Notice that for
each $\tau\in [0,T]\setminus M$ we can write $h_v'(\tau,v)=(
g'_v(\tau,v+u(\tau),\eps)-g'_v(\tau,v_0,0))-(g'_v(\tau,v,0)-g'_v(\tau,v_0,0))$.
As a direct consequence of (\ref{newdef}) we deduce that
$\|h_v'(\tau,v)\|\leq 2\widetilde{\gamma}$ for all $v\in
B_\delta(v_0)$ and $\tau\in [0,T]\setminus M$. Now applying the
mean value theorem for the function $h(\tau,\cdot)$, we have
$\|h(\tau,v_1)-h(\tau,v_2)\|\leq 2\widetilde{\gamma} \|v_1-v_2\|$
for all $\tau\in[0,T]\setminus M$ and all $v_1,v_2\in
B_\delta(v_0)$. Then $\|{F}_2(v_1)-{F}_2(v_2)\|\leq
\int\limits_{[0,T]\setminus M}\|h(\tau,v_1)-h(\tau,v_2)\|d\tau\le
2T\widetilde{\gamma}\|v_1-v_2\|\le(\gamma/2)\|v_1-v_2\|.$
Therefore, we have proved that $\|F(v_1)-F(v_2)\|\leq
\gamma\|v_1-v_2\|,$ that coincides with (\ref{newunif}).
\end{proof}

\begin{lemma}\label{lemma2}
Let $g_0:\Omega\to\R^k$ satisfying assumption (\ref{NE0o}) for some
$v_0\in\Omega.$ Then there exist  $q\in[0,1),$ $\alpha,\, \delta_0>0$ and a norm $\|\cdot\|_0$ on $\R^k$ such that (\ref{conde}) is satisfied.
\end{lemma}

\begin{proof}
If $\lambda$ is an eigenvalue of $\alpha ({g_0})'(v_0)$ then
$\lambda+1$ is an eigenvalue of $I+(\alpha{g_0})'(v_0).$ Since the
eigenvalues of $\alpha({g_0})'(v_0)$ tends to $0$ as $\alpha\to 0$
and have negative real parts then  there exists $\alpha\in[0,1)$
such that the absolute values of all the eigenvalues of
$I+\alpha({g_0})'(v_0)$ are less than one. Therefore (see
\cite[p.~90, Lemma~2.2]{krapo}) there exist $\widetilde{q}\in[0,1)$
and a norm $\|\cdot\|_0$ on $\R^k$ such that
$\sup_{\|\xi\|_0\le 1}\|\xi+\alpha({g}_0)'(\ups_0)\xi\|_0\le\widetilde{q}.$

By continuous differentiability of $g_0$ in a neighborhood of $v_0$
we have that $\| g_0(v_1)-g_0(v_2)-(g_0)'(v_0)(v_1-v_2)\| \, /\, \|
v_1-v_2 \| \le
\|g_0(v_1)-g_0(v_2)-(g_0)'(v_2)(v_1-v_2)\|+\|(g_0)'(v_2)(v_1-v_2)-(g_0)'(v_0)(v_1-v_2)\|/\|v_1-v_2\|
\to 0$ as $\max\{\|v_1-v_0\|,\|v_2-v_0\|\} \to 0.$ Therefore taking
into account that all norms on $\R^k$ are equivalent, there exists
$\delta_0>0$ such that $ \|
g_0(v_1)-g_0(v_2)-(g_0)'(v_0)(v_1-v_2)\|_0 \leq
(1-\widetilde{q})/(2\alpha)\, \| v_1-v_2 \|_0$ for all $v_1,v_2\in
B_{\delta_0}(v_0)$. Then
$$
\begin{array}{rl}
\|v_1+&\hskip-0.3cm\alpha{g}_0(v_1)-v_2- \alpha{g}_0(v_2)\|_0\\
\leq & \alpha \| g_0(v_1)-g_0(v_2)-(g_0)'(v_0)(v_1-v_2)\|_0 +
\|v_1-v_2+\alpha({g}_0)'(v_0)(v_1-v_2)\|_0
\\
\leq & (1+\widetilde{q})/2 \, \|v_1-v_2 \|_0,
\end{array}
$$
for all $v_1,v_2\in B_{\delta_0}(v_0)$.
\end{proof}

{\it Proof of Theorem~\ref{th1}.} Lemmas~\ref{lemma1} and
\ref{lemma2} imply that assumptions (\ref{newunif}) and
(\ref{conde})  of Theorem~\ref{th0} are satisfied. Therefore the
conclusion of the theorem follows applying Theorem~\ref{th0}. $\Box$

It was observed by Mitropol'skii in \cite{mit} that in spite of the
fact that $g(t,\cdot,\eps)$ in (\ref{ps}) is only Lipschitz,
function $g_0$ turns out to be differentiable in applications. In
particular, one will see in Section~3 that this is the case for the
nonsmooth van der Pol oscillator.

Clearly if $g\in
C^1(\mathbb{R}\times\mathbb{R}^k\times[0,1],\mathbb{R}^k)$ then
(\ref{lip}) and (\ref{newdef}) hold in any open bounded set
$\Omega\subset\mathbb{R}^k$.
Therefore Theorem~\ref{th1} is a generalization of the periodic
case of the second Bogolyubov's theorem formulated in the
introduction.

\begin{remark}\label{ad2}
Our Theorem~\ref{th1} does not require that the eigenvectors of
$(g_0)'(v_0)$ be orthogonal as in the result of Buic\u{a} and
Daniilidis  (\cite{adr}, Theorem~3.6). Moreover assumption
($H_2$) of \cite{adr} is more restrictive than (\ref{newdef}).
\end{remark}

For completeness we give also the following theorem on the
existence of non--asymptotically stable $T$--periodic solutions for
(\ref{ps}). In the theorem below, $d(F,V)$ denotes the Brouwer
topological degree of the vector field $F\in
C^0(\mathbb{R}^k,\mathbb{R}^k)$ on the open and bounded set
$V\subset\R^k$ (see \cite[Ch.~2, \S~5.2]{kraop}).

\begin{theorem}\label{th15}
Let $g\in C^0(\mathbb{R}\times\mathbb{R}^k\times[0,1],\mathbb{R}^k).$ Assume
that there exists an open bounded set $V\subset\R^k$ such that
$g_0(v)\not=0$ for any $v\in\partial V$ and
\begin{enumerate}
\setcounter{enumi}{\value{items}}
\item\label{neg}
$d(-{g}_0,V)<0.$ \setcounter{items}{\value{enumi}}
\end{enumerate}
Then there exists $\eps_0>0$ such that for any $\eps\in(0,\eps_0]$
system (\ref{ps}) has at least one non--asymptotically stable
$T$--periodic solutions $x_\eps$ with $x_\eps(0)\in V.$
\end{theorem}

\begin{proof} Since $g_0(v)\not=0$ for any $v\in\partial V$
then from Mawhin's Theorem \cite{maw1} (or \cite[Section~5]{maw3})
we have that there exists $\eps_0>0$ such that
\begin{equation}\label{co}
  d(-{g}_0,V)=d(I-x(T,\cdot,\eps),V)\quad{\rm for\ any\
  }\eps\in(0,\eps_0].
\end{equation}
By \cite[Theorem~9.6]{kraop} for any asymptotically stable
$T$--periodic solution $x_\eps$ of (\ref{ps}) we have that
$d(I-x(T,\cdot,\eps),B_\delta(x_\eps(0)))=1$ for $\delta>0$ sufficiently small. Therefore if all the possible
$T$--periodic solutions of (\ref{ps}) with $\eps\in(0,\eps_0]$ had
been  asymptotically stable, then the degree
$d(I-x(T,\cdot,\eps),V)$ would have been nonnegative,
contradicting (\ref{neg}) and (\ref{co}).
\end{proof}

\begin{remark}
Assumptions (\ref{def}) and (\ref{conde}) imply that
$d(-{g}_0,V)=1$ (see \cite[Theorem~5.16]{kraop}).
\end{remark}

Finally thinking in the application to the nonsmooth van der Pol
oscillator, we formulate the following theorem which combines
Mawhin's Theorem (see \cite{maw1} or \cite[Theorem~3]{maw3},
Theorem~\ref{th1} and Theorem~\ref{th15}. In this theorem
$([g_0]_i)'_{(j)}$ stays for the derivative of the $i$--th component
of the function $g_0$ with respect to the $j$--th variable.

\begin{theorem}\label{th2}  Let $g\in
C^0(\mathbb{R}\times\Omega\times[0,1],\mathbb{R}^2).$  Let
$v_0\in\Omega$ be such a point that ${g}_0(v_0)=0$ and $g_0$ is
continuously differentiable in a neighborhood of $v_0.$
\begin{itemize}
\item[(a)] If ${\rm det}\,({g}_0)'(v_0)\neq 0$
then there exists $\eps_0>0$ such that for any $\eps\in(0,\eps_0]$
system (\ref{ps}) has at least one $T$--periodic solution $x_\eps$
such that $x_\eps(0)\to v_0$ as $\eps\to 0.$
\item[(b)] If
(\ref{lip}) and (\ref{newdef}) hold
 and
\begin{equation}\label{DE1}
{\rm det}\,({g}_0)'(v_0)>0 \quad \mbox{and} \quad
{\left([g_0]_1\right)}'_{(1)}(v_0)+{\left(
[g_0]_2\right)}'_{(2)}(v_0)<0,
\end{equation}
then there exists $\eps_0>0$ such that for any $\eps\in(0,\eps_0]$
system (\ref{ps}) has exactly one $T$--periodic solution
$x_\eps$ such that $x_\eps(0)\to v_0$ as $\eps\to 0.$ Moreover
the solution $x_\eps$ is asymptotically stable.
\item[(c)] If ${\rm det}\,({g}_0)'(v_0)<0$,
then there exists $\eps_0>0$ such that for any $\eps\in(0,\eps_0]$
system (\ref{ps}) has at least one non--asymptotically stable
$T$--periodic solution $x_\eps$ such that $x_\eps(0)\to v_0$ as
$\eps\to 0.$
\end{itemize}
\end{theorem}

\begin{proof} Statement (a) is added for the completeness of the
formulation of Theorem~\ref{th2} and  it follows from Mawhin's
Theorem (see \cite{maw1} or \cite[Theorem~3]{maw3}).

On the other hand it is a simple calculation to show that
(\ref{DE1}) implies that all the eigenvalues of $({g}_0)'(v_0)$
have negative real part. Therefore, assumption (\ref{NE0o}) of
Theorem~\ref{th1} is also satisfied and  statement (b) follows
from this theorem.

Statement (c) follows from Theorem~\ref{th15}. Indeed since ${\rm det}\,({g}_0)'(v_0)<0$ implies (see \cite[Theorem~5.9]{kraop}) that
$d({g}_0,B_\rho(v_0))$ is defined for any $\rho>0$ sufficiently
small and that $d({g}_0,B_\rho(v_0))={\rm det}({g}_0)'(v_0)<0.$
\end{proof}

\section{Application to the nonsmooth van der Pol oscillator}

In his paper \cite{hog} Hogan first demonstrated the existence of
a limit cycle for the nonsmooth van der Pol equation $\ddot
u+\eps(|u|-1)\dot u+u=0$ which governs the circuit drawn at
Fig.~\ref{thom} with the lamp characteristic $i(u)=S_0+S_1
u-S_2v|v|$ whose derivative $i'(u)=S_1-2S_2|v|$ is
nondifferentiable (see Nayfeh-Mook \cite{nay}, \S3.3.4,  where the
same stiffness characteristic appears in mechanics). In this paper
we extend this study considering the van der Pol problem on the
location of stable and unstable periodic solutions of the
perturbed equation
\begin{equation}\label{vp}
\ddot u+\eps \left(|u|-1\right)\dot u+(1+a\eps)u=\eps\lambda\sin t,
\end{equation}
where $a$ is a detuning parameter and $\eps\lambda\sin t$ is an
external force. We discuss with respect to the parameters $a$ and
$\lambda$, under the assumption that $\eps>0$ is sufficiently
small.

Levinson's change of variables (see \cite{lev}, passing from
Eq~2.0 to Eq~2.1) allows to rewrite equation (\ref{vp}) in a
smooth form where the second Bogolubov's theorem is applied. But
we remind that the idea of considering this example is to see what
is the issues of the direct applying of Theorem~\ref{th1} in
comparison with the smooth results by Andronov and Vitt.

Some function $u$ is a solution of (\ref{vp}) if and only if
$(z_1,z_2)=(u, \dot u)$ is a solution of the system
\begin{equation}\label{psvp}
\begin{array}{lll}
    \dot z_1 & = & \dot z_2, \\
    \dot z_2 & = & -z_1+\eps[-az_1-(|z_1|-1)z_2+\lambda\sin t].
    \end{array}
\end{equation}
After the change of variables
$$
  \left(\begin{array}{cc}
    z_1(t) \\
    z_2(t)
    \end{array}\right)=\left(\begin{array}{cc}
    \cos t & \sin t \\
    -\sin t & \cos t
    \end{array}\right)\left(\begin{array}{cc}
    x_1(t) \\
    x_2(t)
    \end{array}\right),
$$
system (\ref{psvp}) takes the form
\begin{equation}\label{psvp1}
  \begin{array}{lll}
  \dot x_1&=&\eps\sin(-t)\left[-a(x_1\cos t+x_2\sin t)-\right.\\ & &\left.-\left(|x_1\cos
  t+x_2\sin t|-1\right)(-x_1\sin t+x_2\cos t)+\lambda\sin t\right],\\
 \dot x_2&=&\eps\cos(-t)\left[-a(x_1\cos t+x_2\sin t)-\right.\\ & &-\left.\left(|x_1\cos
  t+x_2\sin t|-1\right)(-x_1\sin t+x_2\cos t)+\lambda\sin t\right].
  \end{array}
\end{equation}
The corresponding
averaging function $g_0$, calculated according to the formula
(\ref{averfunction}), is given by
\begin{equation}\label{axeq}
  \begin{array}{lll}
    [g_0]_1(M,N)&=&\pi a N-\pi\lambda+\pi
    M-\frac{4}{3}M\sqrt{M^2+N^2},\\
    {[g_0]}_2(M,N)&=&-\pi a M+\pi
    N-\frac{4}{3}N\sqrt{M^2+N^2},
  \end{array}
\end{equation}
and it is continuously differentiable in
$\mathbb{R}^2\backslash\{0\}$.

\begin{figure}[h]
\begin{center}
\includegraphics[scale=0.63]{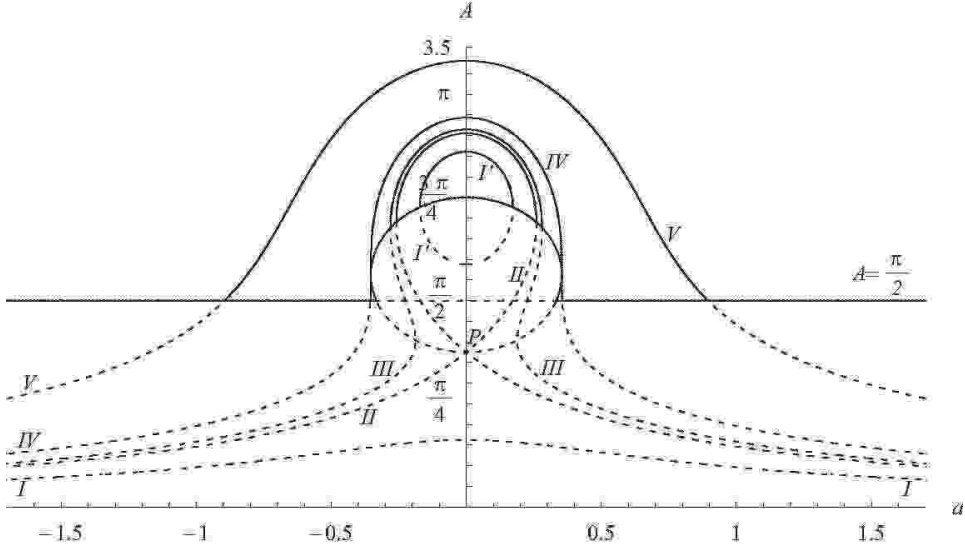}\end{center}
\caption{\footnotesize Dependence of the amplitude of stable
(solid curves) and unstable (dash curves) $2\pi$--periodic
solutions of the nonsmooth periodically perturbed van der Pol
equation (\ref{vp})  on the detuning parameter $a$ obtained over
formulas (\ref{eqA}), (\ref{A1}) and (\ref{A2}) for different
values of $\lambda.$  The curve $I$ is plotted with $\lambda=0.4,$
$II$ with $\lambda=3\pi/16,$ $III$ with some
$\lambda=\sqrt{0.4}\in(3\pi/16,9\sqrt{3}\pi/64),$ $IV$ with
$\lambda=9\sqrt{3}\pi/64,$ $V$ with $\lambda=1.5.$ Point $P$ is
$2/\sqrt{3}.$} \label{fig1}
\end{figure}

In short, by statement~(a) of Theorem~\ref{th2}, the zeros
$(M,N)\in \R^2$ of this function with the property that ${\rm
det}\,(g_0)'(M,N)\not=0$,
 determine the $2\pi$--periodic solutions of
(\ref{psvp}) emanating from the solution of the unperturbed system
\begin{equation}\label{su}
  \begin{array}{lll}
    u_1(t)&=&M\cos t+N\sin t,\\
    u_2(t)&=&-M\sin t+N\cos t.
  \end{array}
\end{equation}
  One has the following expression for the
determinant
\begin{equation}\label{de}
{\rm
det}\,(g_0)'(M,N)=\pi^2(1+a^2)+\frac{32}{9}(M^2+N^2)-4\pi\sqrt{M^2+N^2}.
\end{equation}
Following Andronov and Witt \cite{andr} we are concerned with the
dependence of the amplitude of the solution (\ref{su}) with
respect to $a$ and $\lambda$, thus we decompose this solution as
follows
\begin{equation}\label{eq}
  u_1(t)=A\sin(t+\phi),\ \ u_2(t)=A\cos(t+\phi),
\end{equation}
where $(M,N)$ is related to $(A,\phi)$  by
\begin{equation}\label{A-P}
M=A\sin\phi,\ \ N=A\cos\phi.
\end{equation}
Substituting (\ref{A-P}) into (\ref{axeq}) and (\ref{de}) we obtain
\begin{equation}\label{r}
\begin{array}{l}
 [g_0((A\sin\phi,A\cos\phi))]_1 = -(4/3)\cdot A|A|\sin\phi+\pi aA\cos\phi+\pi
  A\sin\phi-\pi\lambda,\\
   \left[g_0((A\sin\phi,A\cos\phi))\right]_2  = -(4/3)\cdot A|A|\cos\phi-\pi aA\sin\phi+\pi
  A\cos\phi
\end{array}
\end{equation}
and, respectively
\begin{equation}\label{de1}
{\rm
det}\,(g_0)'((A\sin\phi,A\cos\phi))=\pi^2(1+a^2)+\frac{32}{9}A^2-2\pi|A|.
\end{equation}
Looking for the zeros $(A,\phi)$ of (\ref{r}), we find the following
implicit formula for determining $A$:
\begin{equation}\label{eqA}
  A^2\left(a^2+\left(1-\frac{4}{3\pi}|A|\right)^2\right)=\lambda^2.
\end{equation}
Observe that the number of positive zeros of equation (\ref{eqA})
coincides with the number of zeros of the equation
$A^2\left(a^2+\left(1-\frac{4}{3\pi}A\right)^2\right)=\lambda^2.$ To estimate this number we define
$$
f(A)=A^2\left(a^2+\left(1-\frac{4}{3\pi}A\right)^2\right)-\lambda^2,
$$
and we have
$$
f'(A)=2A\left(a^2+\left(1-\frac{4}{3\pi}A\right)^2\right)-
\frac{8}{3\pi}A^2\left(1-\frac{4}{3\pi}A\right).
$$
Since $f'$ has one or two zeros then equation (\ref{eqA}) has one,
two or three positive solutions $A$ for any fixed $a$ and
$\lambda.$ In order to understand the different situations that
can appear, we follow Andronov and Witt who suggested in
\cite{andr} (see also \cite{andr1}) to construct the so called
{\it resonance curves}, namely the curves of dependence of $A$ on
$a$, for fixed $\lambda$. Formula (\ref{eqA}) is the equation of
this curve. Some curves (\ref{eqA}) corresponding to different
values of $\lambda$ are drawn in Figure~\ref{fig1}. The way for
describing these resonance curves (\ref{eqA}) is borrowed from
\cite[Ch.~1, \S~5]{mal}, where the classical van der Pol equation
is considered.

When $\lambda=0$ the curve (\ref{eqA}) is formed by the axis $A=0$
and the isolated point $(0,3\pi/4).$ When $\lambda>0$  but
sufficiently small the resonance curve consists of two branches:
instead of $A=0$ we have the curve of the type $I-I$ and instead of
the point $(0,3\pi/4)$ we obtain an oval $I'-I'$ surrounding this
point. When $\lambda>0$ increases, the oval $I'-I'$ and the branch
$I-I$ tend to each other and, for a certain $\lambda$ there exists
only one branch $II-II$ with a double point $P.$ The value of this
$\lambda$ can be obtained assuming that equation (\ref{eqA}) has for
$a=0$ a double root and, therefore, (\ref{de1}) should be zero.
Solving jointly (\ref{eqA}) and (\ref{de1}) with $a=0$ we obtain
$\lambda=3\pi/16$ and $P=2\pi/8.$ If $\lambda>3\pi/16$ then we have
curves of the type $III$ which take form $V$ when $\lambda>0$
crosses the value $\lambda=9\sqrt{3}\pi/64.$ From here, if
$\lambda<3\pi/16$, then equation (\ref{eqA}) has three real roots
when $|a|$ is sufficiently small, and only one root when $|a|$ is
greater than a certain number which depends on $\lambda.$ When $3\pi/16<\lambda<9\sqrt{3}\pi/64$ we will have one, three
or one solution according to whether $a<a_1,$ $a_1<a<a_2$ or
$a>a_2,$ where $a_1,a_2$ depend on $\lambda.$ The amplitude curves
of  type $V$ provide exactly one solution of (\ref{eqA}) for any
value of $a.$ The value $\lambda=9\sqrt{3}\pi/64,$ that separates
the curves where (\ref{eqA}) has three solutions from the curves
where (\ref{eqA}) has one solution, is obtained from the property
that (\ref{eqA}) with this $\lambda$ has a double root for some $a$
and thus this value of $a$ vanishes (\ref{de1}). Therefore
$\lambda=9\sqrt{3}\pi/64$ is the point separating the interval
$(0,\lambda)$ where the system formed by (\ref{eqA}) and
\begin{equation}\label{de0}
  \pi^2(1+a^2)+\frac{32}{9}A^2-2\pi|A|=0
\end{equation}
has at least one solution from the interval $(\lambda,\infty)$
where (\ref{eqA})--(\ref{de0}) has no solutions.

In short we have studied the amplitudes of the $2\pi$--periodic
solutions of  system (\ref{psvp}) depending on $a$ and $\lambda.$
Whether a physical system described by (\ref{psvp}) possesses
$2\pi$--periodic oscillations corresponding to some of these
$2\pi$--periodic solutions depend on whether some of these
$2\pi$--periodic solutions are asymptotically stable. To find the
answer we use statement~(b) of Theorem~\ref{th2}. Assumption
(\ref{lip}) is obviously satisfied with $\Omega=\mathbb{R}^2.$  Next
statement shows that the right hand side of system (\ref{psvp1})
satisfies (\ref{newdef}).

\begin{proposition}\label{newprop}
Let $v_0\in \mathbb{R}^2$, $v_0\neq 0$. Then the right hand side of
(\ref{psvp1}) satisfies (\ref{newdef}) for any $a,\lambda\in\mathbb{R}.$
\end{proposition}

The proof of the proposition is given in an appendix after this
section.

Thus we have to study the signs of (\ref{de1}) and
${\left([g_0]_1\right)}'_M(A\sin\phi,A\cos\phi)\linebreak+
{\left([g_0]_2\right)}'_N(A\sin\phi,A\cos\phi).$
We have
$$
{\left([g_0]_1\right)}'_M(M,N)+{\left([g_0]_2\right)}'_N(M,N)=
2\left(\pi-2\sqrt{M^2+N^2}\right),
$$
and therefore the conditions for the asymptotic stability of the
$2\pi$--periodic solutions of (\ref{psvp}) near (\ref{su}) are
\begin{equation}\label{depos}
\pi^2(1+a^2)+\frac{32}{9}(M^2+N^2)-4\pi\sqrt{M^2+N^2}>0,
\end{equation}
and
\begin{equation}\label{muneg}
  2\left(\pi-2\sqrt{M^2+N^2}\right)<0.
\end{equation}
Substituting (\ref{A-P}) into the inequalities (\ref{depos}) and
(\ref{muneg}), we obtain the following equivalent inequalities in
terms of the amplitude $A$
\begin{equation}\label{A1}
\pi^2(1+a^2)+\frac{32}{9}A^2-2\pi|A|>0
\end{equation}
and
\begin{equation}\label{A2}
  2\pi-4|A|<0.
\end{equation}
Conditions (\ref{A1}) and (\ref{A2}) mean that the asymptotically
stable $2\pi$--periodic solutions of (\ref{psvp}) correspond to
those parts of resonance curves under consideration which are
outside of the ellipse (\ref{de0}) and above the line $A=\pi/2.$
All the results are collected in Figure~\ref{fig1} from where it
is easy to see that for any detuning parameter $a$ and any
amplitude $\lambda>0$,  equation (\ref{vp}) possesses at least one
asymptotically stable $2\pi$--periodic solution with amplitude
close to $A$ obtained from (\ref{eqA}). Among all the
asymptotically stable $2\pi$--periodic solutions of (\ref{vp}),
there exists exactly one whose fixed neighborhood does not contain
any non--asymptotically stable $2\pi$--periodic solution of
(\ref{vp}) for sufficiently small $\eps>0.$ The amplitude of this
asymptotically stable $2\pi$--periodic solution is obtained from
(\ref{A1})--(\ref{A2}).

\begin{figure}[h]
\begin{center}\includegraphics[scale=0.63]{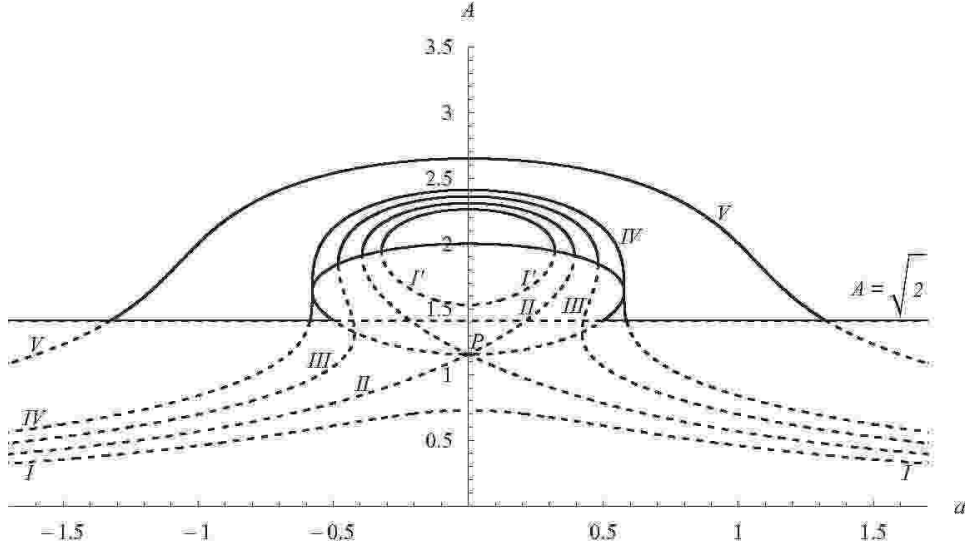}\end{center}
\caption{\footnotesize Dependence of the amplitude of stable
(solid curves) and unstable (dash curves) $2\pi$--periodic
solutions of the classical periodically perturbed van der Pol
equation (\ref{vpor}) on the detuning parameter $a$ for different
values of $\lambda.$ Following Andronov--Witt (see \cite{andr},
Fig.~4) curve $I$ is plotted with $\lambda=\sqrt{0.4},$ $II$ with
$\lambda=4\sqrt{3}/9,$ $III$ with some
$4\sqrt{3}/9<\lambda<\sqrt{32/27},$ $IV$ with
$\lambda=\sqrt{32/27},$ $V$ with $\lambda=2.$ Point $P$ is
$2/\sqrt{3}.$} \label{fig2}
\end{figure}

To compare the changes due to nonsmoothness in the behavior of the
resonance curves,
 we give in Figure~\ref{fig2}
 the resonance curves of the
classical van der Pol oscillator
\begin{equation}\label{vpor}
  \ddot u+\eps \left(u^2-1\right)\dot
  u+(1+a\eps)u=\eps\lambda\sin t,
\end{equation}
which can be found in \cite[Fig.~4]{andr} or in \cite[Ch.~I, \S~16,
Fig.~15]{mal}.

\noindent The key formulas for Figure~\ref{fig1} can be also
easily comparable with the key formulas for Figure~\ref{fig2}. In
fact the corresponding expressions (\ref{eqA})--(\ref{de0}) and
(\ref{depos})--(\ref{muneg}) are (see the formulas~(5.21)--(5.22)
and (16.6)--(16.7) from \cite{mal})
\[
A^2\left(a^2+\left(1-\frac{A^2}{4}\right)^2\right)=\lambda^2,
\]
\[
1-a^2-A^2+\frac{3}{16}A^4=0,
\]
and
\[
1+a^2-(M^2+N^2)+\frac{3}{16}(M^2+N^2)^2>0,
\]
\[
2-(M^2+N^2)<0,
\]
respectively, when considering the classical van der Pol equation
(\ref{vpor}).

\smallskip

It can be checked that the eigenvectors of the
matrix $(g_0)'((A\sin\phi,A\cos\phi))$ are orthogonal only for $A=0$
that is Theorem~3.6 from Buic\u{a}--Daniilidis paper \cite{adr} can
not be applied. At the same time,  assumption~(H2) from \cite{adr} is not satisfied for our problem (see Remark~\ref{ad2}).

\section{Appendix}

{\it Proof of Proposition~\ref{newprop}.} As before, $[v]_i$ is
the $i$--th component of the vector $v\in\mathbb{R}^2.$ Let
$g(t,v)=\left|[v]_1\cos t+[v]_2\sin t\right|$ and notice that it
is enough to prove that $g:[0,2\pi]\times \R^2 \to \R$ satisfies
(\ref{newdef}). In the case that $[v_0]_2\neq 0$, denote
$\theta(v)=\arctan (-[v]_1/[v]_2)$, while when $[v_0]_2=0$, denote
$\theta(v)=\arctan (-[v]_1/[v]_2)$ for $[v_0]_1[v]_2<0$,
$\theta(v_0)=\pi/2$ and, respectively, $\theta(v)=\arctan
(-[v]_1/[v]_2)+\pi$ for $[v_0]_1[v]_2>0$. In any case notice that
the function $v\mapsto \theta(v)$ is continuous in every
sufficiently small neighborhood of $v_0$. Fix
$\widetilde{\gamma}>0.$ Let $M$ be the union of two intervals
centered in $\theta(v_0)$ (when $\theta(v_0)<0$, take
$\theta(v_0)+2\pi$ instead) and, respectively, $\theta(v_0)+\pi$,
each of length $\widetilde{\gamma}/2$. Denote them $M_1$ and
$M_2$. Take $\widetilde{\delta}>0$ such that $\theta(v)\in M_1$
for all $v\in B_{ \widetilde{\delta}}(v_0)$. Of course, also
$\theta(v)+\pi\in M_2$ for all $\|v-v_0\|\leq \widetilde{\delta}$.
This implies that for fixed $t\in [0,2\pi]\setminus M$, $[v]_1
\cos t +[v]_2 \sin t$ has constant sign for all $v\in B_{
\widetilde{\delta}}(v_0)$, that, further, gives that $g(t,\cdot)$
is differentiable and $ g'_v(t,v)=g'_v(t,v_0)$ for all $v\in B_{
\widetilde{\delta}}(v_0)$. Hence (\ref{newdef}) is fulfilled.
$\Box$


\smallskip

\section*{Acknowledgements}
The second author is partially supported by a \linebreak
MEC/FEDER grant number MTM2005-06098-C02-01 and by a CICYT grant
number 2005SGR 00550 and the third author is partially supported
by the Grant BF6M10 of Russian Federation Ministry of Education
and CRDF (US), and by RFBR Grants 06-01-72552, 05-01-00100.

A part of this work was done during a visit of the first and the
third author at the Centre de Recerca Matem\`{a}tica, Barcelona
(CRM). They express their gratitude to the CRM for providing very
nice working conditions.

The authors are grateful to the anonymous referee who motivated to
include a list of several applications of Theorem~\ref{th1} in the
introduction, that definitely improved the paper. The authors also
thank Aris Daniilidis for helpful discussions  and Rafael Ortega
who brought to their attention the change of variables suggested
in the Levinson's paper \cite{lev}.

Finally, we thank M. Golubitsky and A. Vanderbauwhede who invited
us to present the paper at their minisimposia "Recent developments
in bifurcation theory" of Equadiff 2007 (see \cite{blm}), that
gave a significant impact to its recognition.

\end{document}